\documentclass[12pt]{amsart}
\usepackage[margin=1 in]{geometry}
\usepackage{subcaption}
\usepackage{amssymb}
\usepackage{prettyref}
\usepackage{amssymb, amsmath,amsfonts,amsthm,euscript,mathrsfs,  MnSymbol,verbatim, enumerate,multirow,bbding, multicol,
color,array,tikz, tikz-cd,tikz-3dplot,ytableau,multicol}
\usetikzlibrary{calc}
\usepackage{hyperref}
\usepackage{tikz}

\usepackage[labelformat=simple]{subcaption}

\newtheorem{theorem}{Theorem}[section]
\newtheorem*{theorem*}{Theorem}
\newtheorem{proposition}[theorem]{Proposition}
\newtheorem{corollary}[theorem]{Corollary}
\newtheorem{lemma}[theorem]{Lemma} 

\theoremstyle{definition}
\newtheorem{definition}[theorem]{Definition}

\theoremstyle{remark}
\newtheorem{remark}[theorem]{Remark}
\newtheorem{example}[theorem]{Example}

\numberwithin{equation}{section}

\newcommand{\poset}{\mathbb{E}}
\newcommand{\posetz}{\mathbb{E}^\ast}
\newcommand{\poseti}{ \widehat{\mathbb{E}}}

\newcommand{\Exterior}{\mathchoice{{\textstyle\bigwedge}}%
    {{\bigwedge}}%
    {{\textstyle\wedge}}%
    {{\scriptstyle\wedge}}}
    
\DeclareMathOperator{\erays}{Rays}

\newcommand{\glngeo}{I(\mathfrak{gl}_n,V \oplus \bigwedge^2)}
\newcommand{\glngeoForTitle}{\texorpdfstring{$\glngeo$}{I(gl(n))} }

\newcommand{\glngeoSm}[1]{I(\mathfrak{gl}_{#1},V \oplus \bigwedge^2)}

\newcommand{\rayFn}{\vec{v}}
\newcommand{\ray}[1]{\rayFn\left(#1 \right)}
\newcommand{\indFn}{\vec{v}^{-1}}
\newcommand{\ind}[1]{\indFn\left(#1 \right)}

\newcommand{\pref}{\prettyref}
\newrefformat{assert}{(\ref{#1})}
\newrefformat{prop}{Proposition \ref{#1}}
\newrefformat{cor}{Corollary \ref{#1}}
\newrefformat{def}{Definition \ref{#1}}
\newrefformat{rem}{Remark \ref{#1}}
    
\begin{document}
\title[Incidence Geometry in a   Weyl Chamber I: $GL_n$ ]
 {Incidence Geometry in a Weyl Chamber I: $GL_n$}
 
\author[Esole]{Mboyo Esole}
\author[Jackson]{Steven Glenn Jackson}
\author[Jagadeesan]{Ravi Jagadeesan}
\author[No\"el]{Alfred G. No\"el}

\address{Department of Mathematics\\
Harvard University\\ Cambridge, MA 02138, USA.}
\email{esole@math.harvard.edu}
\email{rjagadeesan@college.harvard.edu}
\address{Department of Mathematics\\
         University of Massachusetts\\
         Boston, MA 02125, USA.}
\email{jackson@math.umb.edu}
\email{anoel@math.umb.edu}

\begin{abstract}
We study the  central hyperplane arrangement whose hyperplanes are  
the vanishing loci of the  weights of the  first and the second fundamental representations of  $\mathfrak{gl}_n$ restricted to the dual  fundamental Weyl chamber. 
We obtain  generating functions  that count  flats and  faces of a given dimension. 
This counting is interpreted in physics as  the enumeration of the  phases of the Coulomb and mixed Coulomb-Higgs branches of a five dimensional gauge theory with 8 supercharges 
in presence of hypermultiplets transforming in the fundamental and antisymmetric 
representation of a $U(n)$ gauge group as  described by the Intriligator-Morrison-Seiberg superpotential.
\end{abstract}
\maketitle

\setcounter{page}{1}
 \tableofcontents{}

\newpage

\section{Introduction and summary of results}\label{sec:intro}

The theory of hyperplane arrangements 
has deep connections with many areas of mathematics, 
such as combinatorics, representation theory, algebraic geometry, algebraic topology, singularity theory, 
and the theory of hypergeometric functions.  One beautiful application of the use of hyperplane arrangements in algebraic geometry is the proof of Brieskorn's conjecture by Deligne~\cite{MR0422673} in 1972. 
Real hyperplane arrangements serve as the motivation for  the axioms that define oriented matroids. See \cite{MR1217488,MR2383131,MR1311028} for detailed reviews of the theory of hyperplane arrangements.

Real hyperplane arrangements have been used recently in the physics literature to analyze the fiber structure and  the network of flop transitions between different resolutions of elliptically fibered Calabi-Yau threefolds \cite{
  ESY1,ESY2,Hayashi:2014kca}. 
Compactifications of M-theory on Calabi-Yau threefolds form an elegant bridge  between the study of supersymmetric gauge theories in five dimensional spacetime and  certain types of highly structured  hyperplane arrangements defined  by a Lie algebra $\mathfrak{g}$ and a representation $\mathbf{R}$ of $\mathfrak{g}$.
These hyperplane arrangements have the peculiarity of being defined  not in the full affine  space as they are 
restricted to  the  dual fundamental Weyl chamber of the Lie algebra $\mathfrak{g}$.

\begin{definition}\label{def:incGeoGeneral}
Let $\mathfrak{g}$ be a reductive Lie algebra over $\mathbb{C}$, let $\mathfrak{h}$ be a split, real form of a Cartan subalgebra of $\mathfrak{g},$ and let ${\bf R}$ be a representation of $\mathfrak{g}$.
We denote by $\mathrm{I}(\mathfrak{g},{\bf R})$ the real hyperplane arrangement consisting of the kernels of the weights of ${\bf R}$ restricted to a dual fundamental Weyl chamber in $\mathfrak{h}$.
\end{definition}

\begin{remark}
Because all Cartan subalgebras are conjugate and all fundamental Weyl chambers canonically related by the Weyl group action, the incidence geometry $\mathrm{I}(\mathfrak{g},{\bf R})$ is independent of the choice of $\mathfrak{h}$ and of a dual fundamental Weyl chamber in $\mathfrak{h}$.
\end{remark}

\begin{remark}
 In the rest of the paper, we call the ``dual  fundamental Weyl chamber" the {\em Weyl chamber}. 
\end{remark}

An arrangement  of hyperplanes stratifies the ambient space into open polyhedral cones called {\em faces}. 
A non-empty  intersection of a finite number of  hyperplanes of the arrangement is called a {\em flat}. A $k$-face  (resp. $k$-flat) is a face  (resp.  a flat) that generates a $k$-dimensional linear space.  Faces of maximum dimension are called \emph{chambers}.
A classical combinatorial problem asks how many $k$-flats and $k$-faces a given hyperplane arrangement contains. 
Zaslasky \cite{MR0357135} enumerated chambers  and bounded chambers of an affine hyperplane arrangement using the M\"obius function of its semi-lattice of flats. 

The purpose of this paper is to count the numbers of $k$-faces and $k$-flats of hyperplane arrangements 
$$
I(\mathfrak{gl}_n,V\oplus \Exterior^2),
$$
where $V$ is the vector representation of $\mathfrak{gl}_n$ (the first fundamental representation)  and $\bigwedge^2$ is its second exterior power (the second fundamental representation). 
Since our arrangements are restricted to the Weyl chamber, we cannot apply Zaslasky's result directly.
However, the highly symmetric nature of the hyperplane arrangements under consideration allows us to study them via their extreme rays.

We also analyze the geometries $I(\mathfrak{g},{\bf R})$ when $\mathfrak{g}$ is not simply laced. 
The case of  $I(\mathfrak{sl}_n,V\oplus \Exterior^2)$ is discussed in  \cite{EJJN2}.

\subsection{Motivation from string geometry}
\label{sec:physicsMotivation}
The counting problem solved in this paper is formulated purely in terms of representation theory. However, its motivation comes from the geometry of string theory. 
In this subsection,  we review the motivation.  The reader uninterested  in the interface of string theory and algebraic geometry can safely   ignore this subsection. 

Hyperplane arrangements of the form  $I(\mathfrak{g},\mathbf{R})$ arise naturally in the study \emph{Coulomb branch} of a gauge theory in five-dimensional supersymmetric theory with eight supercharges \cite{Ferrara:1996wv,Intriligator:1997pq}.
Such theories contain different types of particles organized into representations of the supersymmetric algebra called {\em supermultiplets}.   We will consider only  vector multiplets and hypermultiplets. When a gauge theory is defined in a five dimensional supersymmetric theory, the vector multiplets transform in the adjoint representation while the hypermultiplets transform in some representation $\mathbf{R}$ of $\mathfrak{g}$. 
The hypermultiplets  can become massive through the Higgs mechanism along hyperplanes  passing through the origin of the Weyl chamber of $\mathfrak{g}$. These hyperplanes partition the Weyl chamber into different connected regions called \emph{Coulomb phases} of the gauge theory, which are characterized by  a cubic prepotential depending on real coordinates  $\varphi^a$ of  the weight space.

Intriligator-Morrision-Seiberg computed the prepotential of a five dimensional gauge theory for a theory in which hypermultiplets transform in the representation $\mathbf R$ of a   Lie  algebra $\mathfrak{g}$. The result is the real function \cite{Intriligator:1997pq}: 
$$
{\mathcal F}(\varphi) = \cdots 
 +\frac{1}{12}\Big(\sum_{\alpha} | \alpha\cdot \varphi|^3 + \sum_{\lambda} |\lambda\cdot \varphi|^3\Big),
$$
depending on the real vector $\varphi$ in the weight space. We only wrote the terms of ${\mathcal F}(\varphi)$ that are relevant for our purposes. The variable  $\alpha$ runs through the simple roots of $\mathfrak{g}$, and $\lambda$ through the weights of the representation ${\bf R}$ under which the hypermultiplets transform.
The cubic absolute values generate singularities, which correspond to hyperplanes along which new massless particles are generated. 

A well-defined phase of the Coulomb branch of the theory is a connected region in which the quantities $\alpha\cdot \varphi$ and $\lambda\cdot \varphi$ each take fixed signs. 
We fix the sign of $\alpha\cdot \varphi$ by requiring $\varphi$ to be in the fundamental Weyl chamber. 
The condition that $\lambda\cdot \varphi\neq 0$ for all $\lambda\in \mathbf{R}$  restricts us  inside a  specific chamber of the arrangement $I(\mathfrak{g}, \mathbf{R})$. 
The counting of flats is connected to the enumeration of mixed \emph{Coulomb-Higgs branches} of the gauge theory.

A  five dimensional gauge  theory with 8 supercharges can be geometrically engineered by  a compactification of M-theory on a Calabi-Yau threefold \cite{Ferrara:1996wv,Intriligator:1997pq}. When the Calabi-Yau threefold is elliptically fibered, 
its singular fibers determine the   Lie algebra   $\mathfrak{g}$ and the matter representation $\mathbf{R}$  of the gauge theory. 
 A smooth elliptic fibration with a rational section is a resolution of  a singular Weierstrass model, and flop transitions join pairs of crepant resolutions of the Weierstrass model.  
  It is conjectured that the 
  network of (partial) crepant resolutions of the Weierstrass model is isomorphic to the adjacency graph of the faces of the hyperplane arrangement $I(\mathfrak{g},\mathbf{R})$.  
  More precisely, it is conjectured that each crepant resolution corresponds to a unique chamber of $I(\mathfrak{g},\mathbf{R})$ \cite{ ESY1,ESY2,Hayashi:2014kca, Intriligator:1997pq}, and   that two crepant resolutions are connected by a flop if and only if the intersection of the corresponding chambers is a face of codimension one. Their common face is then conjectured to correspond to a unique  partial resolution of the Weierstrass model. See \cite{ESY1,ESY2,EY1} for examples of  explicit matchings between (partial) resolutions and faces of a hyperplane arrangements $I(\mathfrak{g},\mathbf{R})$. There are subtleties in presence of non-trivial $\mathbb{Q}$-factorial terminal singularities, see for example \cite{EKR1}. 
   A proof of this conjecture will likely rely heavily on ideas from Mori's program along the lines of \cite{MR1259932}.
  The Intriligator-Morrision-Seiberg prepotential  $\mathscr{F}(\varphi)$  can also be obtained geometrically \cite{Diaconescu:1998cn,ES1,Intriligator:1997pq}.

In the case of $\glngeo$, the gauge group is $U(n)$, and therefore the relevant elliptic fibrations each have a Mordell-Well group of rank one  and a singular fiber of type I$^s_n$ over a divisor of the discriminant locus. 
  Models for elliptic fibrations with a Mordell-Weil group of rank one are discussed in \cite{Morrison:2012ei} and \cite{Esole:2014dea}.  Elliptic fibrations for geometric engineering of  $U(n)$ gauge theories have been extensively studied recently in the physics literature    (see  for example  \cite{Braun:2013nqa,Cvetic:2013nia, Kuntzler:2014ila,Lawrie:2015hia,Mayrhofer:2012zy,Morrison:2014era} and references therein).

\subsection{Enumeration of faces and flats in \glngeoForTitle}
We are now ready to discuss the enumerations performed in this paper.

We denote  the number of $k$-faces of $I(\mathfrak{gl}_n,V\oplus \Exterior^2)$  by  $g(n,k)$ and collect these numbers in the  following formal power series (generating function):
\[
G(s,t)=\sum_{n,k\geq 0} g(n,k) s^n t^k.\]
The following theorem enumerates faces of $\glngeo$. It is  proved in \pref{sec:genFuncFaces}.

\begin{theorem}
\label{thm:faceCount} The generating function counting faces of $I(\mathfrak{gl}_n, V\oplus \bigwedge{}^2)$ by dimension is given by the following rational function:
\[
G(s,t) =\sum_{n,k\geq 0} g(n,k) s^n t^k= \frac{1 - s}{1 - 2 s + s^2 - 2 s t + s^2 t}.
\]
\end{theorem}

\pref{thm:faceCount} implies a simple recurrence relation for $g(n,k)$.

\begin{corollary}
\label{cor:faceRecurse}
The sequence $g(n,k)$ obeys the recurrence relation:
$$g(n,k) = 2g(n-1,k)-g(n-2,k)+2g(n-1,k-1)-g(n-2,k-1),$$
for $n\geq 2$.
\end{corollary}

We can also obtain formulas for the generating polynomials
\[
G_n(t) = [s^n] G(s,t) = \sum_{k=0}^n g(n,k)t^k.
\]

\begin{corollary}
\label{cor:faceGenPoly}
The generating polynomials $G_n(t)$ are given by
\[
G_n(t)=\frac{\left(t+\sqrt{t (1+t)}\right) \left(1+t+\sqrt{t (1+t)}\right)^n-\left(t-\sqrt{t (1+t)}\right) \left(1+t-\sqrt{t (1+t)}\right)^n}{2 \sqrt{t(1+t)}}.
\]
\end{corollary}

\begin{example}
\pref{cor:faceGenPoly} yields the following values of $G_n(t)$:
\vskip .3 cm
\noindent
\scalebox{.9}{
$
\begin{aligned}
G_0(t) &=1 \\
G_1(t) &=1+t \\
G_2(t) &=1+ 5t+4 t^2\\
G_3(t) &=1+9 t+16 t^2 +8 t^3\\ 
G_4(t) &=1+14 t +41 t^2+44 t^3+16 t^4\\
G_5(t) &=1+20 t+85 t^2 +146 t^3+112 t^4+32 t^5\\
G_6(t) &=1+27 t + 155 t^2+377 t^3+456 t^4+272 t^5+64 t^6\\
G_7(t) &=1+35 t+259 t^2+833 t^3+1408 t^4+1312 t^5+640 t^6+128 t^7\\
G_8(t) &= 1+44 t+406 t^2+1652 t^3+3649 t^4+4712 t^5+3568 t^6+1472 t^7+256 t^8\\
G_{9}(t) &= 1+54 t+606 t^2+3024 t^3+8361 t^4+14002 t^5+14608 t^6+9312 t^7+3328 t^8+512 t^9\\
G_{10}(t) &= 
1+65 t+870 t^2+5202 t^3+17469 t^4+36365 t^5+48940 t^6+42800 t^7+23552 t^8+7424 t^9
 +1024 t^{10}
\end{aligned}
$
}

From the above formulae, we can read off the value of $g(n,k)$ for $k \le n \le 10$ as the coefficient of $t^k$ in $G_n(t)$:
\[
\begin{array}{c|c|c|c|c|c|c|c|c|c|c|c}
g(n,k) & 0 & 1 & 2 & 3 & 4 & 5 & 6 & 7 & 8 & 9 & 10\\ \hline
0 & 1\\ \hline
1 & 1 & 1\\ \hline
2 & 1 & 5 & 4 \\ \hline
3 & 1 & 9 & 16 & 8\\ \hline
4 & 1 & 14 & 41 & 44 & 16 \\ \hline
5 & 1 & 20 & 85 & 146 & 112 & 32 \\ \hline
6 & 1 & 27 & 155 & 377 & 456 & 272 & 64 \\ \hline
7 & 1 & 35 & 259 & 833 & 1408 & 1312 & 640 & 128 \\ \hline
8 & 1 & 44 & 406 & 1652 & 3649 & 4712 & 3568 & 1472 & 256 \\ \hline
9 & 1 & 54 & 606 & 3024 & 8361 & 14002 & 14608 & 9312 & 3328 & 512 \\ \hline
10 & 1 & 65 & 870 & 5202 & 17469 & 36365 & 48940 & 42800 & 23662 & 7424 & 1024
\end{array}.
\]
\end{example}

We can also extract formulae for $g(n,k)$ for $k$ small and $n-k$ small from \pref{thm:faceCount}.

\begin{corollary}
\label{cor:faceSmallLarge}
We have
\[
g(n,k) = \sum_{i=0}^k (-1)^{k-i}2^i\binom{k}{i}\binom{n+i}{2k}.
\]
and
\[
g(n,n-k) = \sum_{i=0}^{\min\left\{k,\left\lfloor\frac{n+1}{2}\right\rfloor\right\}} (-1)^i2^{n-2i} \binom{n-i}{k-i} \left(\binom{n-i+1}{i} + \binom{n-i}{i-1}\right).\]
\end{corollary}

\begin{example}
\pref{cor:faceSmallLarge} yields the following formulae for $g(n,k)$ and $g(n,n-k)$ for $0 \le k \le 4$:
\begin{align*}
 g(n,0) &=1\\
 g(n,1) &=\frac{1}{3}n(n+3) \\
g(n,2) &=\frac{1}{4!}n(n - 1) (n^2 + 11 n + 22)\\
g(n,3) &=\frac{1}{6!}  n (n-1)(n-2)(n + 7)  (n^2 + 17 n + 36)\\
g(n,4) &=\frac{1}{8!}n(n-1)(n-2)(n-3) (n^4+42 n^3+563 n^2+2754 n+3912)\\
g(n,n-4) &=2^{n-11} (27 n^4-234 n^3 + 697 n^2 -810 n +280)\\
 g(n,n-3) &= 2^{n-7} (9n^3-42 n^2+57 n-20)\\
g(n-2) &= 2^{n-5} ( 9 n^2-17 n + 6)\\
 g(n,n-1) &= 2^{n-2} (3n-1)\\
 g(n,n) &=2^n.
\end{align*}
\end{example}

 \pref{thm:faceCount} and Corollaries ~\ref{cor:faceRecurse},~\ref{cor:faceGenPoly} , and~\ref{cor:faceSmallLarge} are proved in \pref{sec:genFuncFaces}.

We denote  the number of $k$-flats of $I(\mathfrak{gl}_n,V\oplus \Exterior^2)$  by  $h(n,k)$ and collect these numbers in the  following formal power series (generating function):
\[
H(s,t)=\sum_{n,k\ge0} h(n,k) s^n t^k.\]
The following theorem enumerates flats  of $\glngeo$.

\begin{theorem}
\label{thm:flatCount} The generating function counting flats of $I(\mathfrak{gl}_n, V\oplus \bigwedge{}^2)$ by dimension is given by the following rational function:
\[
H(s,t) = \frac{ (1-s) (1-s t +s^2 t)}{1-2 s+s^2-2 s t+3 s^2 t-2 s^3 t+s^2 t^2-2 s^3 t^2+s^4 t^2}.
\]
\end{theorem}
  \pref{thm:flatCount} is proven   \pref{sec:genFuncFlats}.
It implies a simple recurrence relation for $h(n,k)$, the analogue of \pref{cor:faceRecurse} for flats.

\begin{corollary}
\label{cor:flatRecurse}
The sequence $h(n,k)$ obeys the recurrence relation
\begin{multline*}
h(n,k) = 2h(n-1,k)-h(n-2,k)+2h(n-1,k-1)-3h(n-2,k-1)\\
+2h(n-3,k-1)-h(n-2,k-2)+3h(n-3,k-2)-h(n-4,k-2)
\end{multline*}
for $n \ge 4$.
\end{corollary}

\begin{example}
The Taylor expansion of \pref{thm:flatCount} to order 10 in $s$ is:
\vspace{.2cm}

\scalebox{.96}{
$
\begin{aligned}
H(s,t) &= 1 \\
& + (1+ t)s\\
& + (1+ 3 t  + t^2)s^2\\
& + (1 +5  t +  6 t^2 + t^3)s^3\\
& + (1 + 8  t + 14 t^2 + 10 t^3 + t^4)s^4\\
& + (1 + 12  t + 29 t^2 + 31 t^3 + 15 t^4 + t^5)s^5\\
& + (1 + 17 t + 54 t^2 + 79 t^3 + 60 t^4 + 21 t^5 + t^6)s^6\\
& + (1 + 23 t + 93 t^2 + 175 t^3 + 183 t^4 + 106 t^5 + 28 t^6 + t^7)s^7\\
& + (1 + 30 t + 151 t^2 + 352 t^3 + 471 t^4 + 380 t^5 + 175 t^6 + 36 t^7 + t^8)s^8\\
& + (1 + 38 t + 234 t^2 + 659 t^3 + 1082 t^4 + 1119 t^5 + 728 t^6 + 274 t^7 + 45 t^8 + t^9)s^9\\
&  + (1 + 47 t + 349 t^2 + 1166 t^3 + 2286 t^4 + 2894 t^5 + 2426 t^6 + 1310 t^7 + 411 t^8 + 55 t^9 + t^{10})s^{10}.
\end{aligned}
$
}
\vspace{.2cm}

From the above formula, we can read off the value of $h(n,k)$ for $k \le n \le 10$ as the coefficient of $s^nt^k$ in $H(s,t)$:
\[
\begin{array}{c|c|c|c|c|c|c|c|c|c|c|c}
h(n,k) & 0 & 1 & 2 & 3 & 4 & 5 & 6 & 7 & 8 & 9 & 10\\ \hline
0 & 1\\ \hline
1 & 1 & 1\\ \hline
2 & 1 & 3 & 1\\ \hline
3 & 1 & 5 & 6 & 1\\ \hline
4 & 1 & 8 & 14 & 10 & 1\\ \hline
5 & 1 & 12 & 29 & 31 & 15 & 1\\ \hline
6 & 1 & 17 & 54 & 79 & 60 & 21 & 1\\ \hline
7 & 1 & 23 & 93 & 175 & 183 & 106 & 28 & 1\\ \hline
8 & 1 & 30 & 151 & 352 & 471 & 380 & 175 & 36 & 1\\ \hline
9 & 1 & 38 & 234 & 659 & 1082 & 1119 & 728 & 274 & 45 & 1\\ \hline
10 & 1 & 47 & 349 & 1166 & 2286 & 2894 & 2426 & 1310 & 411 & 55 & 1
\end{array}.
\]
\end{example}

The form of the series $H(s,t)$ does not allow simple explicit formulae for $H_n$ or for $h(n,k)$ as in Corollaries~ \ref{cor:faceGenPoly} and~\ref{cor:faceSmallLarge}, respectively.  \pref{thm:flatCount} is proved in \pref{sec:genFuncFlats}.

\subsection{Structure of chambers, faces, flats, and extreme rays of \glngeoForTitle}

\begin{figure}
\ytableausetup{nosmalltableaux, boxsize=2.5em}
\begin{ytableau}
(1,1)  & (1,2) & (1,3) & \cdots & (1,n)\\
 \none &(2,2)& (2,3) & \cdots & (2,n) \\
  \none &\none & (3,3) & \cdots & (3,n)\\
  \none & \none & \none & \ddots & \vdots \\
  \none & \none & \none & \none & (n,n)
  \end{ytableau}
\caption{{\bf Sign tableau notation.} The box labeled by $(i,j)$ corresponds to  the sign of  weight $x_i+x_j$.  This convenient  
   parameterization of chambers of $\glngeo$  was introduced in the physics literature \protect\cite{Hayashi:2014kca}.}
\label{fig:signTableau}
\end{figure}

A face is  determined by assigning signs to $x_i+x_j$ for all values of $(i,j)$ such that  $i,j\in \{1,2,\cdots, n\}$. 
These weights can be organized in a  right-justified Young tableau as in Figure \ref{fig:signTableau}.
It is important to realize that not all sign patterns are allowed.
Since we are in the Weyl chamber, we have
\[
i<j\iff x_i-x_j>0.
\]
Using the identities
\begin{align*}
x_i+x_{j+1} &=x_i+x_j-(x_{j}-x_{j+1})\\
x_{i+1}+x_{j} &=x_i+x_j- (x_{i}-x_{i+1}),\\
\end{align*}
one can prove the following simple sign rules 
 called {\em sign flows} in the physics literature  \cite{Hayashi:2014kca}:
\begin{proposition}[\protect{\cite[Section 2.4]{Hayashi:2014kca}}]
The signs in a tableau corresponding to a chamber satisfy the following conditions.
\begin{itemize}
\item All the boxes above or on the left of a box with positive entries are also positive. 
\item All the boxes below or on the right of a box with negative  entries are also  negative. 
\end{itemize}
\end{proposition}

A tableau satisfying these two rules corresponds to a unique chamber  
of $I(\mathfrak{gl}_n, V\oplus \Exterior^2)$.
However, a different notation for chambers will be more convenient for dealing with extreme rays.

\begin{remark}
A different notation for chambers will be more convenient for dealing with extreme rays.
The sign rules implies that the positive entries of each row are next to each other and start on the left border of the table.  It is therefore efficient to denote a given chamber by the numbers of positive entries on each row. 
The sign rules is then automatically satisfied if these numbers form a decreasing sequence 
$\{a_1>a_2> \cdots> a_k\}$, 
where we denote by $a_i$ the number of  positive entries on the $i$th row of the tableau.  
We don't count the rows that do not have any positive entries. 
An entry of the tableau located on the $i$th row and the $j$th column  (with $i\leq j$) is positive if and only if $i\leq k$ and $j\leq a_i$. It is negative otherwise. This justifies the following definition. 
\end{remark}

\begin{definition}
\label{def:subsetToChamber}
For $S = \{a_1 > \cdots > a_k\} \subseteq [n]$, define a face $C(S)$ of $I(\mathfrak{gl}_n, V\oplus \bigwedge{}^2)$ as the subset of $\mathfrak{h}$ on which
\begin{align*}
x_i + x_j \ge 0 &\qquad \text{if } i \le k \text{ and } j \le a_i\\
x_i + x_j \le 0 &\qquad \text{otherwise}.
\end{align*}
We can equivalently write a subset $S \subseteq [n]$ as a \emph{characteristic vector} $s \in \{-1,1\}^n$, where $s = (s_1,\ldots,s_n)$ and $s_i = 2\chi_S(n) - 1$. 
  A characteristic vector $s = (s_1,\ldots,s_n)$ defines a subset $S\subset [n]$ such that  $i\in S$  (resp. $i\notin S$) if  and only if  $s_i=1$  (resp. $s_i=-1$). 
\end{definition}

The notation of \pref{def:subsetToChamber} allows us to state the following classification of chambers and extreme rays.

\begin{theorem}
\label{thm:chamberStructure}
The chambers and extreme rays of $\glngeo$ satisfy the following properties.
\begin{enumerate}[(a)]
\item \label{assert:chamberBiject} The map $S \mapsto C(S)$ defines a bijection from $2^{[n]}$ to the set of chambers of $I(\mathfrak{gl}_n, V\oplus \bigwedge{}^2)$.  In particular, $I(\mathfrak{gl}_n, V\oplus \bigwedge{}^2)$ has $2^n$ chambers.
\item \label{assert:chamberRaysPrelim} The extreme rays of $C(S)$ are generated by the vectors $e_1^S,\ldots,e_n^S,$
 where
\[e_\ell^S = (\underbrace{1,\ldots,1}_{\pi_\ell^S},\underbrace{0,\ldots,0}_{n-\ell},\underbrace{-1,\ldots,-1}_{\ell-\pi_\ell^S}),\]
 where $\pi_\ell^S = \left|S \cap [n-\ell+1,\infty)\right|$   where $\pi_\ell^S = \left|S \cap [n-\ell+1,\infty)\right|$  counts the elements of $S$ that are greater or equal to $n-\ell+1$, and the vectors $e_1^S,\ldots,e_n^S$ non-negatively span $C(S)$. 
In particular, the geometry $\glngeo$ is simplicial.
\end{enumerate}
\end{theorem}

The enumeration of chambers of $\glngeo$ announced in \pref{thm:chamberStructure}\pref{assert:chamberBiject}
recovers the count of phases of $U(n)$ with fundamental and antisymmetric matter from the physics literature \cite{Hayashi:2014kca}.
 Theorem~\ref{thm:chamberStructure} is proven in  \pref{sec:extremeRaysChambers}. 
\begin{remark}
Suppose that $S = \{n_1,\ldots,n_k\}$ with $n_1>n_2>\ldots>n_k$.  A direct calculation shows that the box diagram associated to $C(S)$ has $n_i$ plusses in the $i$th row for $i \le k$ and all minus signs in the $i$th row for $i > k$.
See Tables~\ref{tab:example.gl2} and~\ref{tab:example.gl3} for examples. 
\end{remark}

\ytableausetup{nosmalltableaux, boxsize=1 em}
\begin{table}
\begin{tabular}{c | l|c|  l}
Sign tableau & Subset of $[n]$ & Characteristic vector  & Interior points  \\
\hline
\rule{0pt}{6ex} 
\scalebox{1}{\begin{ytableau}+ & + \\ \none & + \\    \end{ytableau}}  & $\{2,1\}$ &  $(\ \ 1,\ \ 1)$&  $(\ \  a_1+a_2,\ \  a_1)$\\
\rule{0pt}{6ex} 
\scalebox{1}{\begin{ytableau}+ & + \\ \none & - \\    \end{ytableau} } & $\{2\}$ &  $(-1,\ \  1)$&   $(\ \  a_1+a_2,-a_1)$\\
\rule{0pt}{6ex} 
\scalebox{1}{\begin{ytableau}+ & - \\ \none & - \\    \end{ytableau} } & $\{1\}$ & $(\ \ 1,-1)$ &  $(\ \ a_1,-a_1-a_2)$\\
\rule{0pt}{6ex} 
\scalebox{1}{\begin{ytableau}- & - \\ \none & - \\    \end{ytableau}}  & $\emptyset$ & $ (-1,-1)$ &  $(-a_1,-a_1-a_2)$
 \end{tabular}
\vspace{0.5cm}
\caption{{\bf Chambers of $\glngeoSm{2}$.}  We show several different notations for the chambers.  In describing the sets of interior points, the variables $a_1,a_2$ denote positive real numbers.
}
\label{tab:example.gl2}
\end{table}

\begin{table}
\begin{center}
\begin{tabular}{l | c|   c  |     l}
Sign tableau & Subset of $[n]$ & Characteristic vector & Interior points of the  chamber  \\
\hline
\rule{0pt}{6ex} 
\begin{ytableau}- & - & - \\ \none & - & -  \\ \none & \none & -    \end{ytableau} &   $\emptyset$ & $ (-1,-1,-1)$  &   $(-a_1, -a_1-a_2,-a_1-a_2-a_3)$\\
\rule{0pt}{6ex} 
\begin{ytableau} +& - & - \\ \none & - & -  \\ \none & \none & -    \end{ytableau} &   $\{1\}$ & $(\ \ 1,-1,-1)$ & $(a_1,-a_1-a_2 ,-a_1-a_2-a_3)$\\
\rule{0pt}{6ex} 
\begin{ytableau}+ & + & - \\ \none & - & -  \\ \none & \none & -    \end{ytableau}  & $\{2\}$ &  $(-1,\ \  1,-1)$ &  $(a_1+a_2,-a_1 ,-a_1-a_2-a_3)$\\
\rule{0pt}{6ex} 
\begin{ytableau}+ & + & + \\ \none & - & -  \\ \none & \none & -    \end{ytableau}   &$\{3\}$ &  $(-1,-1,\ \ 1)$  &  $(a_1+a_2+a_3,-a_1-a_2 ,-a_1-a_2-a_3)$\\
\rule{0pt}{6ex} 
\begin{ytableau}+ & + & + \\ \none & + & -  \\ \none & \none & -    \end{ytableau}   &$\{3,1\}$ & $(\ \  1,-1,\ \ 1)$ &   $(a_1+a_2+a_3,a_1 ,-a_1-a_2)$ \\
\rule{0pt}{6ex} 
\begin{ytableau}+ & + & - \\ \none & +& -  \\ \none & \none & -    \end{ytableau}   &$\{2,1\}$ &  $(\ \ 1,\ \ 1,-1)$ &  $(a_1+a_2,a_1 ,-a_1-a_2-a_3)$\\
\rule{0pt}{6ex} 
\begin{ytableau}+ & + & + \\ \none & + & +  \\ \none & \none & -    \end{ytableau}   &$\{3,2\}$  &  $(-1,\ \  1,\ \ 1)$ &  $(a_1+a_2+a_3,a_1+a_2 ,-a_1)$ \\
\rule{0pt}{6ex} 
\begin{ytableau}+& + & + \\ \none & + & +  \\ \none & \none & +    \end{ytableau}  & $\{1,2,3\}$ & $(\ \ 1,\  \  1,\ \  1)$ &  $(a_1+a_2+a_3,a_1+a_2 ,a_1)$
 \end{tabular}
\vspace{.5cm}
\end{center}
\caption{{\bf Chambers of $\glngeoSm{3}$.}  We show several different notations for the chambers.  In describing the sets of interior points, the variables $a_1,a_2,a_3$ denote positive real numbers.
} 
\label{tab:example.gl3}
\end{table}

In light of \pref{thm:chamberStructure}, we can introduce additional structure on the set of extreme rays of $\glngeo$.
We will equip the set of extreme rays with a partial order and use the combinatorics of the resulting partially ordered set (poset) to study faces and flats in $\glngeo$.

The following definitions describe the posets that will be relevant to us.

\begin{definition}
\label{def:quarterPlanePoset}
The \emph{discrete quarter plane poset} is the set  $\mathbb{N}^2$ endowed with  the Cartesian order induced by the usual order of the set $\mathbb{N}$ of  non-negative integers: 
$$
 \quad (a,b)\leq (c,d)\iff (a\leq c) \quad and \quad (b\leq d), \quad a,b,c,d\in \mathbb{N}.
$$ 
This is a graded poset with grading function  
$\ell(a,b)=a+b$.  We call $\ell(a,b)$ the \emph{ level}  of $(a,b)$.
\end{definition}

\begin{definition}
\label{def:PosetForExtRays}
We denote by $\poset_n$ the subset of the discrete quarter plane
 that consists of points at  level less or equal to $n$. 
$\poset_n$ is a poset with the order induced  by the Cartesian order defined above. 
We denote by $\poseti_n$ the poset $\poset_n\sqcup \{\infty\}$ where $\infty$ is greater than all the elements of $\poset_n$.    
We denote by $\posetz_n$ the poset $\poset_n$ with the origin removed: 
$$
\begin{aligned}
\poset_n &= \{ (a,b)\in \mathbb{N}^2 \quad | \quad 0\leq a + b\leq n \},\quad
 \quad
\poseti_n  = \poset_n \sqcup \{\infty\}, \quad \posetz_n =\poset_n \setminus \{ (0,0) \}.
\end{aligned}
$$
\end{definition}

We are now ready to relate $\posetz_n$ to the set of extreme rays of $\glngeo$.

\begin{definition}\label{def:vFn}
Define a function $\vec{v}:\posetz_n \to W$ by
\[\ray{a,b}=(\underbrace{1,\cdots ,1}_{a}, \underbrace{0,\cdots ,0}_{n-a-b} , \underbrace{-1,\cdots ,-1}_{b}).\]
\end{definition}

The following theorem follows directly from \pref{thm:chamberStructure}.

\begin{theorem}
\label{thm:rayStruture}
The extreme rays of $I(\mathfrak{gl}_n, V\oplus \bigwedge{}^2)$ satisfy the following properties.
\begin{enumerate}[(a)]
\item The function $\vec{v}$ defines a bijection from $\posetz_n$ to the set of extreme rays of $I(\mathfrak{gl}_n, V\oplus \bigwedge{}^2)$.
\item The set of extreme rays that lie in a chamber $C(S)$ are the extreme rays $\ray{\pi_\ell^S,\ell-\pi_\ell^S},$ where $\pi_\ell^S = \left|S \cap [n-\ell+1,\infty) \right|$. 
\end{enumerate}
\end{theorem}

We will characterize faces and flats by  the extreme rays they contain. For that reason we define a function which returns the extreme rays lying in a given subset of the Weyl chamber.

\begin{definition}
\label{def:raysOperator}
Define a function $\erays$ from $2^W$ to the power set of the set of extreme rays of $\glngeo$ as follows.  For $S \subseteq W,$ let $\erays(S)$ be the set of extreme rays of $\glngeo$ that lie in $W$.
\end{definition}

The following theorem,  proven in \pref{sec:extremeRaysFaces} relates faces in $\glngeo$ to the combinatorics of $\posetz_n$.

\begin{theorem}
\label{thm:Characterization.Faces}
For all $k$, the function $\vec{v}^{-1} \circ \erays$ induces a bijection from the set of $k$-faces of $\glngeo$ to the set of $k$-chains in $\posetz_n$.
\end{theorem}
\begin{remark}
Here, we say that $S \subseteq \posetz_n$ is a $k$-chain if $S$ is a chain and $|S| = k$.
\end{remark}
  
To state an analogue of \pref{thm:Characterization.Faces} for flats, we will need to define structures called \emph{ensembles} that will play the role of chains.
  
\begin{definition}
\label{def:Ensemble}
An \emph{ensemble} is the restriction to $\posetz_n$ of  a union
\[
\left(\bigcup_{0\leq i\leq k}  [A_i,B_i] \right)\cap \posetz_n
\]
of intervals $[A_i,B_i]$ of $\poseti_n$ satisfying the following four conditions:
\begin{enumerate}[(1)]
\item $A_0=0;$
\item $A_i\leq B_i$ for $0 \le i \le k$;
\item $B_i +(1,1) \leq A_{i+1}$ for $0 \le i \le k-1$; and
\item $\ell(B_k)< n$ or $B_k=\infty$.  
\end{enumerate}
We say that $E$ is a \emph{$k$-ensemble} if $k = \left|\ell(E)\right|$, so that $k$ counts the number of distinct levels of elements of $E$.
\end{definition}

\begin{figure}
\begin{multicols}{2}
\begin{minipage}{\linewidth}
\centering
                 \includegraphics[scale=.8]{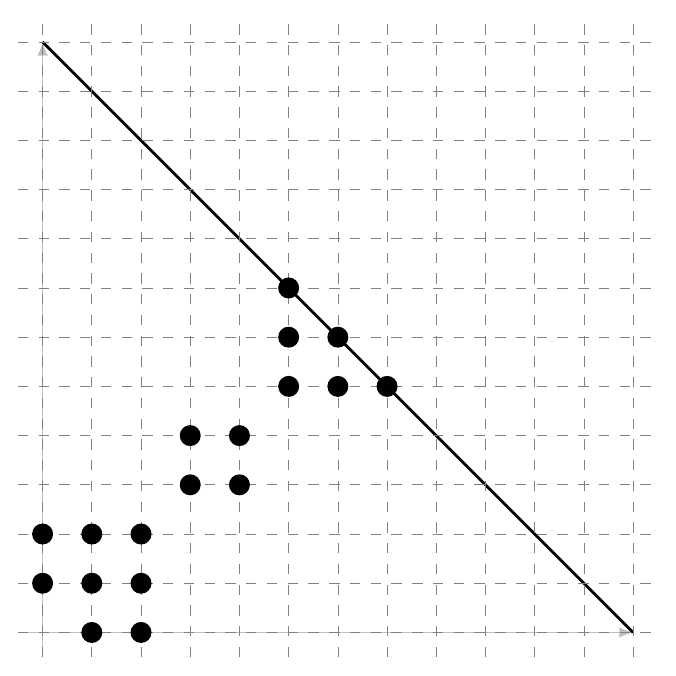}
               \subcaption{  The $10$-ensemble $S$ 
                 in $\posetz_{12}$. \label{fig:fig2A}}
 \end{minipage} 
 
 \begin{minipage}{\linewidth}
\centering
                 \includegraphics[scale=.8]{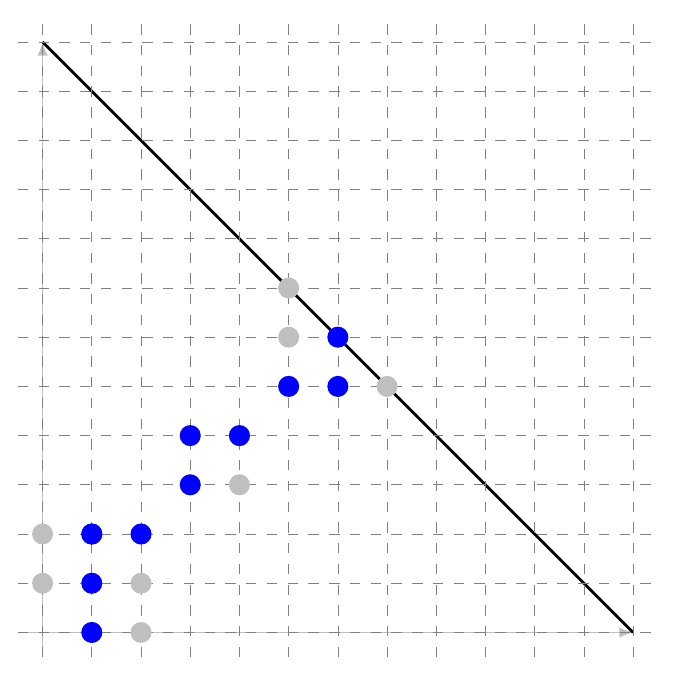}
                 
               \subcaption{A maximal chain of the ensemble $S$. \label{fig:fig2B}}
 \end{minipage} 
 \end{multicols}
\caption{{\bf A $10$-ensemble in $\posetz_{12}$.} 
Figure (A)  represents  the $10$-ensemble $S$ which is  the restriction to $\poset_{12}$ of the union of intervals $[(0,0),(2,2)] \cup [(3,3),(4,4)] \cup [(5,5),\infty]$.
Theorem \ref{thm:Characterization.Flats} guarantees that $\vec{S}$ is the set of extreme rays of a unique 10-flat in $\glngeoSm{12}$.
Figure (B)  illustrates a  maximal chain of  $S$, which corresponds to a $10$-face of $\glngeoSm{12}$. 
}
\label{fig:ensembleOneEg}
\end{figure}

Figure \ref{fig:ensembleOneEg}.\subref{fig:fig2A} shows one example of an ensemble.

\begin{theorem}
\label{thm:Characterization.Flats}
For all $k$, the function $\vec{v}^{-1} \circ \erays$ induces a bijection from the set of $k$-flats of $\glngeo$ to the set of $k$-ensembles of $\posetz_n$.
\end{theorem}
 Theorem \ref{thm:Characterization.Flats} is proven in  section \ref{sec:extremeRaysFlats} on page \pageref{sec:extremeRaysFlats}.
\subsection{Case of a non-simply laced Lie algebra}

Denote by $V$ the vector representation of $\mathfrak{so}_{2n+1}$ and of $\mathfrak{sp}_n$.  Recall that the representation $\bigwedge^2 V$ of $\mathfrak{sp}_n$ factors as $\bigwedge^2_0 \oplus \text{triv}$, where $\bigwedge^2_0$ is irreducible.
Let ${\bf 26}$ (resp. ${\bf 7}$) denote the unique 26-dimensional (resp. 7-dimensional) irreducible representation of $\mathfrak{f}_4$ (resp. $\mathfrak{g}_2$).
We are now ready to discuss the geometry $I\left(\mathfrak{g},{\bf R}\right)$ when $\mathfrak{g}$ is non-simply laced simple Lie algebra.

\begin{theorem}
\label{thm:nonSimplyLaced}
Suppose that
\[(\mathfrak{g},{\bf R}) \in \left\{\left(\mathfrak{so}_{2n+1},V\right),\left(\mathfrak{sp}_n,V\right),\left(\mathfrak{sp}_n,\Exterior^2_0\right),\left(\mathfrak{sp}_n,V \oplus \Exterior^2_0\right),\left(\mathfrak{f}_4,{\bf 26}\right),\left(\mathfrak{g}_2,{\bf 7}\right)\right\}\]
and suppose that $\mathfrak{g}$ has rank $n$.  The geometry $I\left(\mathfrak{g},{\bf R}\right)$ is a cone over a simplex $\Delta^{n-1}$.
In particular, there is 1 chamber and the number of $k$-faces (resp. $k$-flats) is $\binom{n}{k}$.
\end{theorem}

\begin{remark}
There is only one other geometry $I(\mathfrak{g},{\bf R})$ when $\mathfrak{g}$ is simple and non-simply laced and ${\bf R}$ is a direct sum of minuscule and quasi-minuscule representation.
This is the geometry $I\left(\mathfrak{so}_{2n+1},{\bf spin}\right).$
\end{remark}
It is mentioned in  \cite[Section 6.3]{Hayashi:2014kca}  that $I\left(\mathfrak{sp}_n,\bigwedge^2_0\right)$ and $I\left(\mathfrak{sp}_n,V \oplus \bigwedge^2_0\right)$ have no chambers. This is in contradiction with the conclusion of    \pref{thm:nonSimplyLaced}.

\begin{proof}[Proof of \pref{thm:nonSimplyLaced}]
A direct calculation shows that all the weights of ${\bf R}$ are scalar multiples of the roots of $\mathfrak{g}$ in the cases discussed in the theorem.
Therefore, the hyperplanes $\lambda^\perp$ for $\lambda$ a weight of ${\bf R}$ coincide with vanishing loci of roots of $\mathfrak{g}$.  The theorem then follows from the fact that the Weyl chamber is simplicial.
\end{proof}

\begin{figure}
                 \includegraphics[scale=1.1]{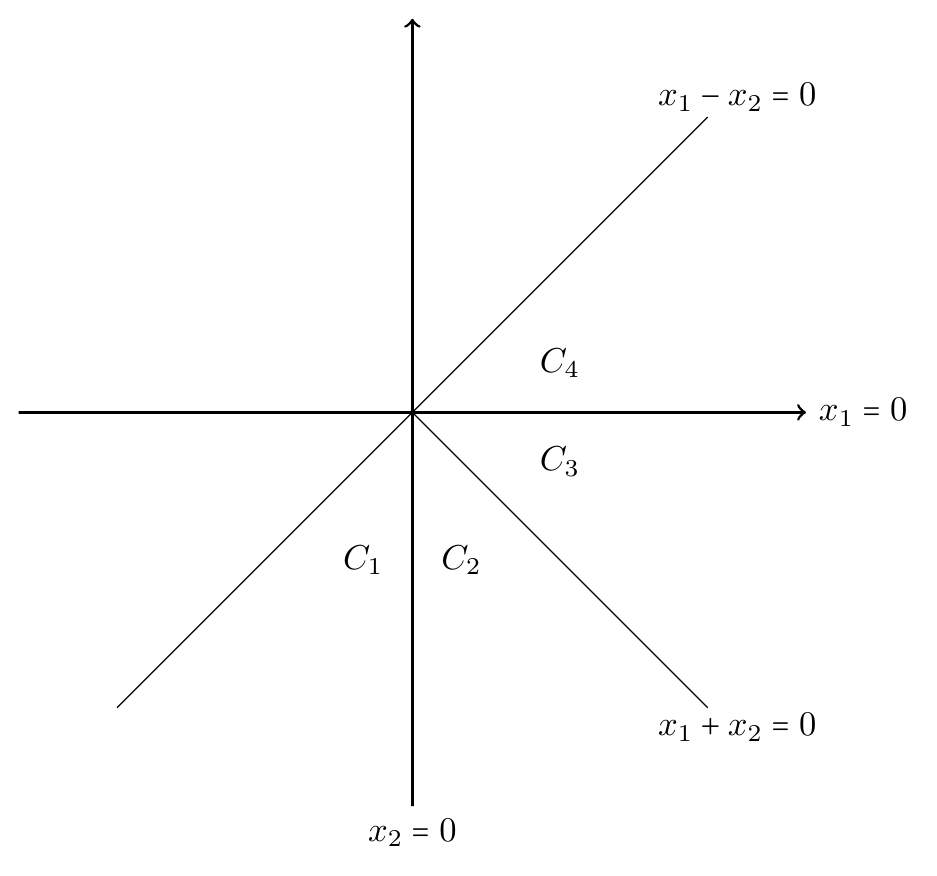}
\caption{
 {\bf The geometry $\glngeoSm{2}$.} The geometry contains  4 chambers, 5 half-lines and the origin. There is one 2-flat (the half-plane below the line $x_1-x_2=0$), three 1-flats (the interior walls $C_1$, $C_2$, $C_3$), and one 0-flat (the origin). }
\label{fig:plane}

\end{figure}

\subsection{Outline of the paper}

In \pref{sec:prelims}, we recall basic facts regarding the Lie algebra $\mathfrak{gl}_n$ and from the theory of hyperplane arrangements.  In \pref{sec:extRaysChmFaces}, we prove Theorems~\ref{thm:chamberStructure},~\ref{thm:rayStruture}, and~\ref{thm:Characterization.Faces}.  In \pref{sec:faceCount}, we prove \pref{thm:faceCount}, and in \pref{sec:flatCount}, we prove Theorems~\ref{thm:flatCount} and~\ref{thm:Characterization.Flats}.

\subsubsection*{Acknowledgments}
M.E. is grateful to Shu-Heng Shao and Shing-Tung Yau for helpful discussions. 
M.E. and A.N. thank William Massey and all the organizers of  the 20th Conference for African American Researchers in Mathematical Sciences (CAARMS 20) where this  started. 
M.E. is supported in part by
the National Science Foundation (NSF) grant DMS-1406925  ``Elliptic Fibrations and String
Theory".
R.J. is supported by the Harvard College Research Program.

\section{Preliminaries}
\label{sec:prelims}
In \pref{sec:prelimLieAlg}, we fix notation for the Weyl chamber of $\mathfrak{gl}_n$.  In \pref{sec:hypArrange}, we recall some facts from the theory of hyperplane arrangements.
A detailed list of symbols is presented in \pref{tab:notations}.

\subsection{The Lie algebra \texorpdfstring{$\mathfrak{gl}_n$}{gl(n)}}
\label{sec:prelimLieAlg}

Let $\mathfrak{h} \subseteq \mathfrak{gl}_n$ denote the Cartan subalgebra of $\mathfrak{gl}_n$ consisting of the diagonal matrices.  We identify $\mathfrak{h}$ with $\mathbb{R}^n$ via the isomorphism
\[(x_1,\ldots,x_n) \leftrightarrow \operatorname{diag}(x_1,\ldots,x_n).\]

The positive roots of $\mathfrak{gl}_n$ are $\{x_i-x_j \mid 1\leq i < j \leq n \}$.
We define the open Weyl chamber of $\mathfrak{gl}_n$ as the set of vectors $v\in \mathfrak{h}$ forming a non-increasing sequence
$$
W^0 =\{ (x_1,\ldots,x_n) \in \mathfrak{h} \mid x_1>\cdots > x_n \}. 
$$
and denote by $W = \overline{W^0}$ the closure of $W^0$.

Define $\lambda_{i,j} = x_i + x_j \in \mathfrak{h}^*$.
The  weights of the vector representation are $x_i$ for   $1 \le i\le n$ and the weights of the antisymmetric representation  
are $x_i+x_j$ for $1 \le i < j \le n$.
When $i=j$, the hyperplane $\lambda^\perp_{i,i}$ is the vanishing locus of the weight $x_i$ of the vector representation of $\mathfrak{gl}_n$.
When $i\neq j$, the hyperplane $\lambda^\perp_{i,j}$ is the vanishing locus of the weight $x_i+x_j$ of the antisymmetric representation of $\mathfrak{gl}_n$.

\subsection{Hyperplane arrangements}
\label{sec:hypArrange}

By definition, each chamber is  a  convex rational polyhedral cone (finite intersection of rational half spaces). In the following we will refer to a  convex rational polyhedral cone simply as ``cone''.
A \emph{face} is either  a chamber or the intersection of a chamber with a supporting hyperplane. There are finitely many faces, and  each face is a lower dimensional cone. 
The chambers cover $W$ and form a \emph{fan}. A fan is a finite collection of  cones such that for each cone all its faces are in the fan and the intersection of two cones is a face of each \cite{MR1234037}.  
 By  classical theorems of Weyl and Minkowski, a cone is finitely generated if and only if it is a finite intersection of closed linear half spaces (see \cite[Theorem 1.3]{MR1311028}). 
It follows  that each chamber of $I(\mathfrak{gl}_n,V\oplus \Exterior^2)$  admits a system of  generators. 
 A cone of dimension $n$  is  said to be \emph{simplical} if it is  non-negatively  spanned by  $n$ linearly independent vectors.  
  A  face of dimension one  which is a half-line is  called an \emph{extreme ray}.  A \emph{facet} is a  face of co-dimension one.

The roots of the Lie algebra  $\mathfrak{gl}_n$  define a hyperplane arrangement known as the {\em braid arrangement} \cite{MR2383131}, which has $n!$ chambers, namely the Weyl chambers of $A_{n-1}$. The braid arrangement is also called the Coxeter arrangement of type $A_{n-1}$. 
The weights of the fundamental representation  define an arrangement of hyperplanes known as the {\em Boolean arrangement} \cite{MR2383131}.
The weights of the antisymmetric  representation $\Lambda^2$ define an arrangement of hyperplanes known as the {\em threshold  arrangement}  \cite{MR2383131}.
Altogether, the roots of $\mathfrak{gl}_n$ and the weights of the first two fundamental representations define the Coxeter arrangement of type  $BC_n$ \cite{MR2133266,MR2360474,MR2383131}. 

\section{Geometry of extreme rays}
\label{sec:extRaysChmFaces}

In \pref{sec:extremeRaysChambers}, we prove Theorems~\ref{thm:chamberStructure} and~\ref{thm:rayStruture}, and in \pref{sec:extremeRaysFaces}, we prove \pref{thm:Characterization.Faces}.

\subsection{Extreme rays and chambers}
\label{sec:extremeRaysChambers}
Given a vector $x\in \mathbb{R}^n$, we denote by $x_{\leq}$ (resp. $x_{\geq}$) the vector obtained from $x$ by  rearranging its components in non-increasing (resp. non-decreasing) order.
The idea of the proof is to derive for each $s \in \{-1,1\}^n$ an explicit linear isomorphism $\alpha_s$ from $\mathbb{R}^n$ to $\mathfrak{h}$ such $\alpha_s(\mathbb{R}_{\ge 0}^n) = C(s)$.
The elementary basis vectors of $\mathbb{R}^n$ then map under $\alpha_s$ to the set of extreme rays of $C(s)$.

\begin{proof}[Proof of \pref{thm:chamberStructure}]
Define
\[\Delta = \bigcup_{i,j} \lambda_{i,j}^\perp.\]
The absolute values of the coordinates  $(x_1, \cdots, x_n)$ of a point $x \in W^0 \setminus \Delta$ are pairwise distinct and non-zero, because $x_i + x_j \not= 0$ and $x_i - x_j \not= 0$ for all $1 \le i < j \le n$. After reordering, we can write  the set $\{|x_1|, |x_2|, \cdots, |x_n|\}$ as  a strictly  increasing vector
\[(|x_1|, |x_2|, \cdots, |x_n|)_{\leq}=(a_1(x), a_1(x) +a_2(x), a_1(x)+a_2(x)+a_3(x), \cdots ,a_1(x) + a_2(x) +\cdots+ a_n(x))\]
so that $a_i(x)\in \mathbb{R}_{>0}$ for all $i$
Because
\[\{a_1(x), a_1(x) +a_2(x), a_1(x)+a_2(x)+a_3(x), \cdots ,a_1(x) + a_2(x) +\cdots+ a_n(x)\} =  \{|x_1|,\ldots,|x_n|\},\]
there  exists a unique sequence $\gamma(x) = (\gamma_1(x),\ldots,\gamma_n(x)) \in \{-1,1\}^n$ of signs such that
\[\gamma_i(a_1(x)+\cdots+a_i(x)) \in \{x_1,\ldots,x_n\}\]
for all $1 \le i \le n$.  We see immediately that  $\gamma$ defines a continuous function from $W^0 \setminus \Delta$ to $\{-1,1\}^n$.
A direct calculation shows that
\[
x=\left(\gamma_1(x)a_1(x), \gamma_2(x) \sum_{i=1}^2 a_i(x), \ldots, \gamma_k(x) \sum_{i=1}^k a_i(x), \ldots,  \gamma_n(x) \sum_{i=1}^n a_i(x)\right)_{\geq}
\]
and that $x \in C(\gamma(x))$.

We now prove Part~\pref{assert:chamberBiject} of the theorem.  Suppose that $s = (s_1,\ldots,s_n) \in \{-1,1\}^n$.  consider the vector
\[v(s) = (s_1, 2s_2, \ldots,ns_n)_{\geq} \in W^0 \setminus \Delta.\]
Note that the absolute values of the components of $v(s)$ are $1,2,\ldots,n$.  We see immediately that $\gamma(v(s)) = s$ for all $s \in \{-1,1\}^n$, so that $v(s) \in C(s)$.  Because $v(s) \in W^0 \setminus \Delta$ and $v(s) \in C(s)$, the face $C(s)$ contains a point of $W^0$ not lying in any hyperplane of $\glngeo$.  It follows that $C(s)$ is a chamber of $\glngeo$ for all $s \in S$.

It remains to prove that $C$ defines a bijection from $\{-1,1\}^n$ to the set of chambers of $\glngeo$.  Consider a chamber $F$ of $\glngeo$ and let $x \in F^0$ be an interior point.  Note that $x \in C(\gamma(x))$.  We see immediately that $F = C(\gamma(x))$ because $C(\gamma(x))$ is a chamber of $\glngeo$.
Suppose that $C(s) = C(t)$.  It follows that $v(s)$ and $v(t)$ lie in the same connected component of $W^0 \setminus \Delta$.  Because $\gamma$ is continuous, we obtain that
\[s = \gamma(v(s)) = \gamma(v(t)) = t,\]
so that $C$ is injective.  Part~\pref{assert:chamberBiject} follows.

We now prove Part~\pref{assert:chamberRaysPrelim}.  Consider a chamber $C(s)$ of $\glngeo$.  Let $S = \{i \mid s_i = 1\}$ denote the subset of $[n]$ corresponding to $x \in S$.
Note that
\begin{align*}
x&=\left(s_1a_1(x), s_2 \sum_{i=1}^2 a_i(x), \ldots, s_k \sum_{i=1}^k a_i(x), \ldots,  s_n \sum_{i=1}^n a_i(x)\right)_{\geq}\\
&=\sum_{\ell=1}^n a_\ell \cdot (\underbrace{0,\ldots,0}_{n-\ell},s_{n+1-\ell},\ldots,s_n)_{\geq}=\sum_{i=1}^n a_\ell e_\ell^{S}
\end{align*}
with $a_i(x) \in \mathbb{R}_{> 0}$.  Hence, $C(s)$ is non-negatively spanned by $\{e_1^S,\ldots,e_n^S\}$. It follows that $e_1^S,\ldots,e_n^S$ are the extreme rays of $C(s)$, as desired.
\end{proof}

\subsection{Extreme rays and faces}
\label{sec:extremeRaysFaces}

Recall the definition of a chain in a poset.

\begin{definition}
A {\em\bf  chain} in a poset is a subset of the poset in which any two elements are comparable.
The \emph{length} of a chain is its size, and a \emph{$k$-chain} is a chain of length $k$.
\end{definition}

It follows from \pref{thm:chamberStructure} that $\indFn \circ \erays$ induces a bijection from the set of chambers of $\glngeo$ to the set of $n$-chains in $\posetz_n$.  To prove \pref{thm:Characterization.Faces}, we will exploit the simpliciality of the chambers of $\glngeo$, which was also proved in \pref{thm:chamberStructure}.

\begin{proof}[Proof of \pref{thm:Characterization.Faces}]
Suppose that $F$ is a $k$-face of $\glngeo$.  Consider a chamber $C$ containing $F$.  Because $C$ is simplicial (by \pref{thm:chamberStructure}), $F$ must be the non-negative span of $k$ extreme rays of $C$.  In particular, we obtain that $\erays(F)$ is a subset of $\erays(C)$ of size $k$.  Because $\ind{\erays(C)}$ is an $n$-chain in $\posetz_n$ and $\ind{\erays(F)}$ is a subset of $\ind{\erays(C)}$ of size $k$, it follows that $\ind{\erays(F)}$ is a $k$-chain in $\posetz_n$.
Therefore, $\indFn \circ \erays$ is a function from the set of $k$-faces of $\glngeo$ to the set of $k$-chains of $\posetz_n$.
The fact that every face is non-negatively spanned by the extreme rays that it contains ensures that $\indFn \circ \erays$ induces an injection from the set of $k$-faces of $\glngeo$ to the set of $k$-chains in $\posetz_n$.

It remains to prove that $\indFn \circ \erays$ is surjective.
Suppose that $E \subseteq \posetz_n$ is a $k$-chain.  Let $F$ denote the non-negative span of $\vec{E}$.  Consider an $n$-chain $\tilde{E}$ containing $E,$ which clearly exists, and the corresponding chamber $C$ (which satisfies the property that $\tilde{E} = \erays(C)$).  It follows from the simpliciality of $C$ (\pref{thm:chamberStructure}) that $F$ is a sub-face of $C$ of dimension $k$ and that $\erays(F) = E$.  Therefore, every $k$-chain in $\posetz_n$ can be expressed in the form $\indFn \circ \erays(F)$ for some $k$-face $F$ of $\glngeo$. The theorem follows.
\end{proof}

\section{Counting faces}
\label{sec:faceCount}

In \pref{sec:PrelimOnPoset}, we describe some basic properties of the posets $\poset_n$ and $\posetz_n$.  In \pref{sec:genFuncFaces}, we prove \pref{thm:faceCount}.

\subsection{Preliminaries on the poset \texorpdfstring{$\posetz_n$}{E}}
\label{sec:PrelimOnPoset}

\begin{figure}
\begin{multicols}{3}
\begin{minipage}{\linewidth}
\centering
                 \includegraphics[scale=1]{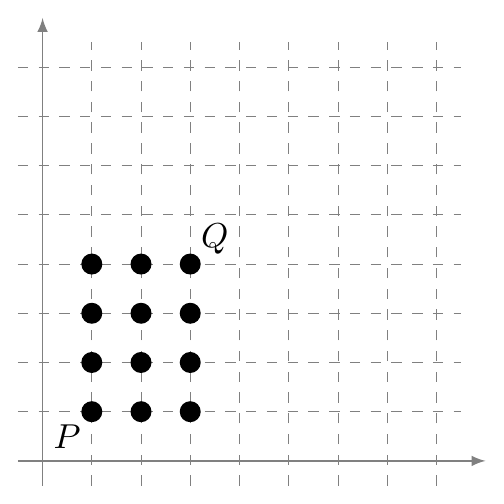}
                 \subcaption{The interval $\lbrack P, Q\rbrack$ is the rectangle suspended between the points $P$ and $Q$. \label{fig:intervalRectangle}}
\end{minipage}

\begin{minipage}{\linewidth}
\centering
                 \includegraphics[scale=1]{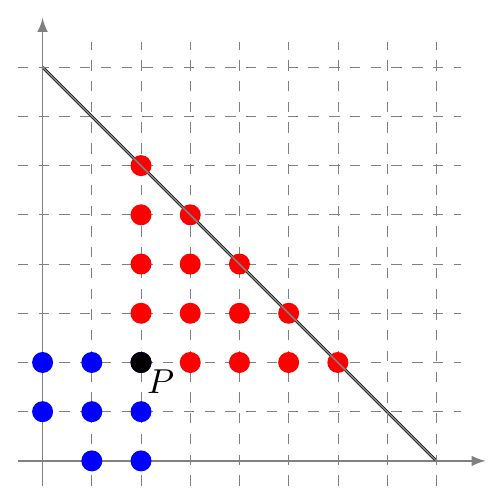}
\subcaption{The blue dots form the lower set  $(O,P)$  of $P$. The red dots form the  upper set $(P,\infty)$ of $P$.\label{fig:intervalLowerUpper}}
\end{minipage}

\begin{minipage}{\linewidth}
\centering
                 \includegraphics[scale=1]{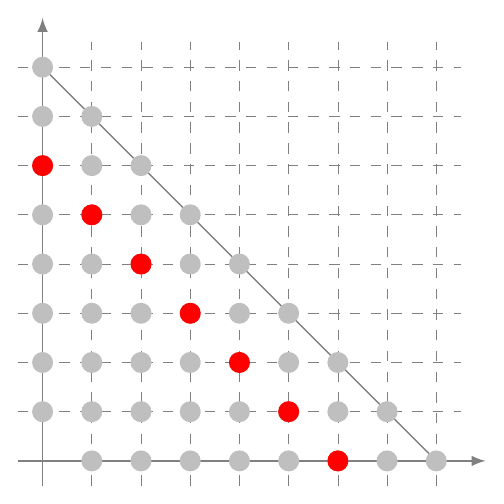}
\subcaption{The red dots form the level set $L_6$, which is a non-empty anti-chain in $\posetz_n$ for $n \ge 6$. \label{fig:levelSet}}
\end{minipage}
\end{multicols}
\caption{{\bf Intervals and level sets in $\poset_n$ and $\posetz_n$.}}
\label{fig:posetsN}
\end{figure}

Given two points $P$ and $Q$ in $\posetz_n$, such that $P\leq Q$,  consider the rectangle with vertical and horizontal edges whose NE-SW  diagonal is the segment $\overline{PQ}$. The  lattice points of this rectangle are exactly the points of the 
interval $[P,Q]$ in $\posetz_n$.  See    Figures  \ref{fig:posetsN}~\subref{fig:intervalRectangle}--\subref{fig:intervalLowerUpper} on page~\pageref{fig:intervalLowerUpper} for examples of intervals in $\posetz_n$.

The discrete quarter-plane $\mathbb{N}^2$ is equipped with a grading given by the level function. The origin is the only point of level zero. 
Denote by $L_i$ the set of points of level $i$ of $\mathbb{N}^2$.
The level sets $L_i$ are all anti-chains and form a partition of $\mathbb{N}^2$.  See \pref{fig:levelSet} for an example. Note that
\[\posetz_n = \bigsqcup_{i = 1}^n L_i \quad \text{ and } \quad \poset_n = \bigsqcup_{i=0}^n L_i.\]
For all $P \in \poset_n$, the interval $(P,\infty)$ in $\poset_n$ (or in $\posetz_n$) is isomorphic to $\posetz_n$ as a graded poset.  The isomorphism is induced by the translation by $P$ in $\mathbb{N}^2$: the map $x \mapsto x + P$ defines an isomorphism from $\posetz_{n-\ell}$ to $(P,\infty)$.

\subsection{Proof of \pref{thm:faceCount}}
\label{sec:genFuncFaces}
In order to compute the generating function $G(s,t),$ we first derive a recurrence relation for $g(n,k)$.

\begin{proposition}
\label{prop:faceRecursion}
The integers $g(n,k)$ satisfy the recurrence relation
\[
g(n,k) = \delta_{k,0}+
\sum_{\ell=1}^{n-k+1} (\ell+1) g(n-\ell,k-1)
\]
for $n \ge 0$, where $g(n,k) = 0$ for $n < 0$ and for $k \notin [0,n]$.
\end{proposition}
We prove \pref{prop:faceRecursion} by tracking the minimum of each $k$-chain and interpreting the $k-1$ largest elements of the chain as $(k-1)$-chain in $\posetz_{n'}$ with $n'<n$.
\begin{proof}
The proposition is clear for $k \le 0$.  Assume for the remainder of the proof that $k \ge 1$.

Denote by $C_P$ the set of $k$-chains in $\posetz_n$ with minimum $P$.
Suppose that $P \in \posetz_n$.  Recall that translation by $-P$ in $\mathbb{N}^2$ defines an isomorphism from $(P,\infty)$ to $ \posetz_{n-\ell(P)}$.  We see immediately that the process of removing $P$ from a $k$-chain and translating by $-P$ (i.e., the function defined by $S \mapsto (S \setminus P) - P$) defines a bijection from $C_P$ to the set of $(k-1)$-chains in $\posetz_{n-\ell(P)}$.  It follows that $\left|C_P\right| = g(n-\ell(P),k-1)$.
Note in particular that $\left|C_P\right| = 0$ for $\ell(p) > n-k+1$.

Recall that $\posetz_n$ contains $\left|L_\ell\right| = \ell+1$ elements of level $\ell$ for $1 \le \ell \le n$.  Grouping by level, we obtain that
\[
g(n,k) =
\sum_{p \in \posetz_n} \left|C_P\right| = 
\sum_{\ell=1}^{n-k+1} (\ell+1)    g(n-\ell,k-1),
\]
as desired.
We need to introduce the correction term $\delta_{k,0}$ in order to deal with the case of $k \le 0$.
\end{proof}

After some algebra, we can derive the expression for $G(s,t)$ announced in \pref{thm:faceCount}.

\begin{proof}[Proof of \pref{thm:faceCount}]
Translating \pref{prop:faceRecursion} in to the language generating functions yields the linear equation
\begin{align*}
G(s,t)=\sum_{ n=0}^\infty \sum_{k=0}^\infty g(n,k) s^n t^k &=\sum_{n=0}^\infty  s^n +\sum_{n=0}^\infty \sum_{k=0}^\infty \sum_{\ell=1}^{\infty} (\ell+1)    g(n-\ell,k-1) s^{n} t^{k} \\
&=\frac{1}{1-s} +t  \Big(\sum_{\ell=1}^{\infty} (\ell+1) s^\ell \Big) \sum_{\ell=0}^\infty \sum_{k=0}^\infty  g(n-\ell,k-1) s^{n-\ell} t^{k-1}\\
&=\frac{1}{1-s} +t \frac{s(2-s)}{(s-1)^2}   \sum_{i=0}^\infty \sum_{j=0}^{\infty}  g(i,j) s^i t^j = \frac{1}{1-s} +t  \frac{s(2-s)}{(s-1)^2}  G(s,t),
\end{align*}
where we used the identity 
$$
\sum_{\ell=1}^{\infty} (\ell+1) s^\ell =\frac{s(2-s)}{(s-1)^2}.
$$
The theorem follows by solving the linear equation for $G(s,t)$.
\end{proof}

Corollaries~\ref{cor:faceRecurse}, \ref{cor:faceGenPoly} and~\ref{cor:faceSmallLarge} follow from \pref{thm:faceCount} by elementary algebraic manipulations.
In particular, \pref{cor:faceRecurse} can be derived by noting that the coefficient of $s^nt^k$ in 
\[(1 - 2 s + s^2 - 2 s t + s^2 t) G(s,t)\]
vanishes for $n \ge 2$ due to \pref{thm:faceCount}.
\pref{cor:faceGenPoly} can be derived by computing a partial fraction decomposition for the expression for $G(s,t)$ given in \pref{thm:faceCount}.  On the other hand, \pref{cor:faceSmallLarge} is proved by regarding $G(s,t)$ as the sum of geometric series with ratio $2 s - s^2 + 2 s t - s^2 t$.

\section{Counting flats}
\label{sec:flatCount}

In \pref{sec:extremeRaysFlats}, we prove \pref{thm:Characterization.Flats}, and in \pref{sec:genFuncFlats}, we prove \pref{thm:flatCount}.

\subsection{Extreme rays and flats}
\label{sec:extremeRaysFlats}
\pref{thm:Characterization.Flats} will be derived from the following two lemmas. The first asserts that $\erays(L)$ is a flat whenever $L$ is an ensemble, and the second asserts that every ensemble is of the form $\erays(L)$ for some flat $L$.

\begin{lemma}
\label{lem:flatYieldsEnsemble}
The set of extreme rays of a $k$-flat of $I(\mathfrak{gl}_n, V\oplus \Exterior^2)$  form a $k$-ensemble. 
\end{lemma}

\begin{lemma}
\label{lem:ensembleYieldsFlat}
Every ensemble of $\poset_n$ is the  set of  extreme rays of a flat of $I(\mathfrak{gl}_n, V\oplus \Exterior^2)$. 
\end{lemma}

\begin{figure}
 
   \begin{multicols}{5}

\begin{minipage}{\linewidth}\centering                 \includegraphics[scale=.7]{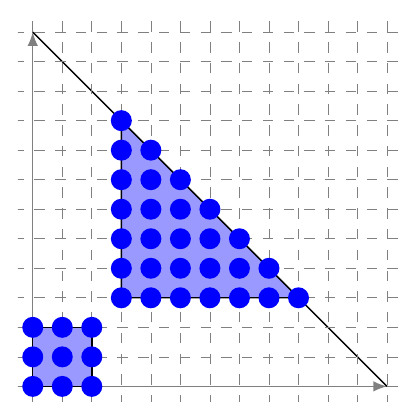}\subcaption{ $\lambda^\bot_{3,n-2}$  \label{fig:fig6a} }\end{minipage}
                 
                 \begin{minipage}{\linewidth}\centering                 \includegraphics[scale=.7]{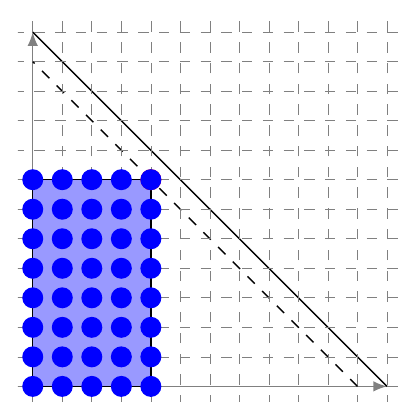}  \subcaption{ $\lambda^\bot_{i,i}$  \label{fig:fig6b} }\end{minipage}

               \begin{minipage}{\linewidth}\centering                 \includegraphics[scale=.7]{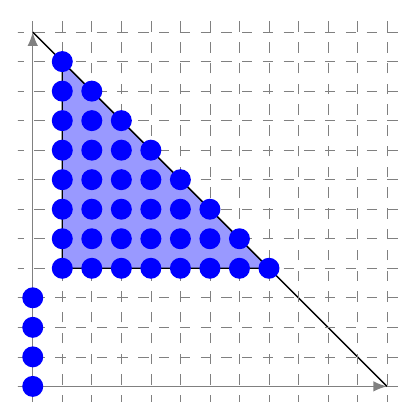}  \subcaption{ $\lambda^\bot_{1,n-3}$  \label{fig:fig6c} }\end{minipage}
               
                                  \begin{minipage}{\linewidth}\centering                 \includegraphics[scale=.7]{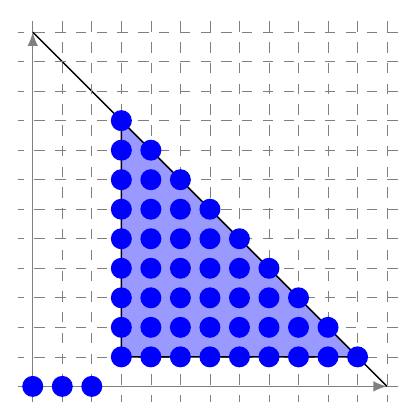}  \subcaption{ $\lambda^\bot_{3,n}$  \label{fig:fig6d} }\end{minipage}
               
                   \begin{minipage}{\linewidth}\centering                 \includegraphics[scale=.7]{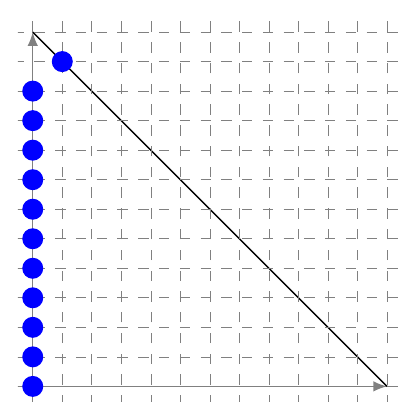}  \subcaption{ $\lambda^\bot_{1,2}$  \label{fig:fig6e} }\end{minipage}
               
                     \begin{minipage}{\linewidth}\centering                 \includegraphics[scale=.7]{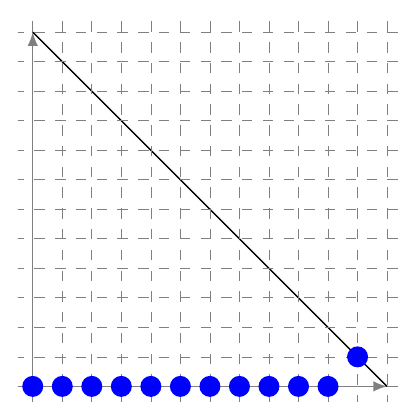}  \subcaption{ $\lambda^\bot_{n-1,n}$  \label{fig:fig6f} }\end{minipage}
               
                       \begin{minipage}{\linewidth}\centering                 \includegraphics[scale=.7]{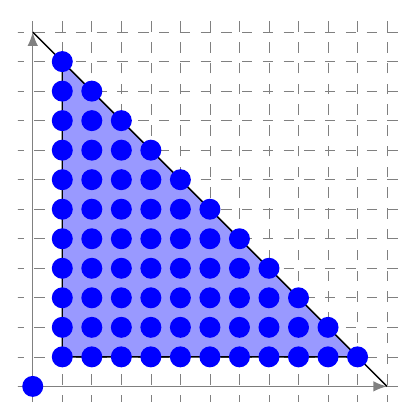}  \subcaption{ $\lambda^\bot_{1,n}$  \label{fig:fig6g} }\end{minipage}
               
                         \begin{minipage}{\linewidth}\centering                 \includegraphics[scale=.7]{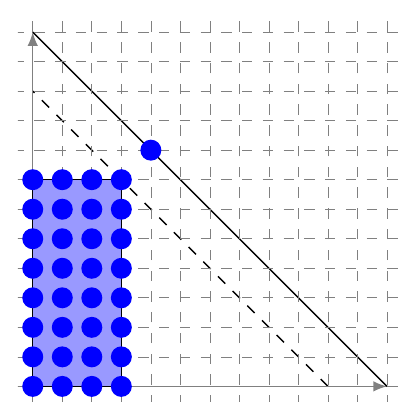}  \subcaption{ $\lambda^\bot_{i,i+1}$  \label{fig:fig6h} }\end{minipage}
              
               \begin{minipage}{\linewidth}\centering                 \includegraphics[scale=.7]{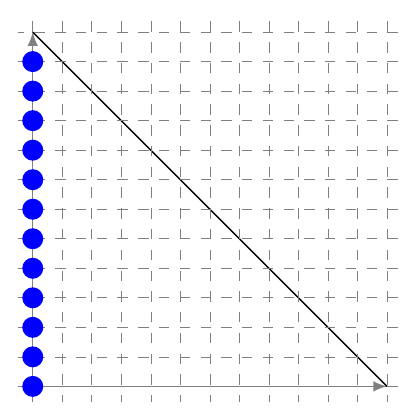}  \subcaption{ $\lambda^\bot_{1,1}$  \label{fig:figi} }\end{minipage}
                              \begin{minipage}{\linewidth}\centering                 \includegraphics[scale=.7]{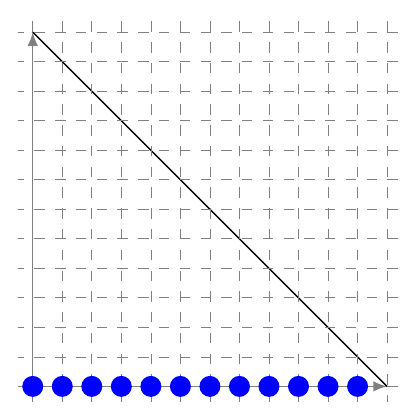}  \subcaption{ $\lambda^\bot_{n,n}$  \label{fig:figj} }\end{minipage}
                \end{multicols}
   
   \caption{ {\bf Set of extreme rays of a hyperplane $\lambda_{i,j}^\bot$}. The set of extreme rays of a hyperplane  $\lambda^\bot_{i,j}$  ($i\leq j$) corresponding to a union of intervals $[0,\varphi(i,j)]\cup [\varphi(i,j)+(1,1),\infty]$ in the poset $\poset_n$ with $\varphi(i,j)=(i-1, n-j)$. 
   For extreme values  ($i=1$, $i=n$, $j=n$, $j=i$, $j=i+1$) one of the two intervals can be contained in one of the axis, shrink to a point, or disappear altogether. 
   The points $(0,n)$ and $(n,0)$ are the only points of $\posetz_n$  not correspond to any extreme ray of any hyperplane $\lambda^\bot_{i,j}$.  
   \label{fig:flats.hyperplane} }
\end{figure}

To prove \pref{lem:flatYieldsEnsemble}, we exploit the the fact that the set of extreme rays of a flat is the intersection of the sets of extreme rays of the hyperplanes containing it.
Recall that the set of extreme rays of $\lambda_{i,j}^\perp$   $(1\leq i\leq j\leq n)$  corresponds in the poset $\posetz_n$ to  the union of intervals
$$(0, \varphi(i,j)]\cup [\varphi(i,j)+(1,1),\infty],$$
where \[\varphi(i,j)=(i-1,n-j).\]
Notice that $\ell(\varphi(i,j)) \le n-1$ with equality if and only if $i = j.$
When $i=j$, the second interval is empty and $\lambda_{i,j}$ is a weight of the fundamental representation $V$. 
The first interval is empty when  $(i,j)=(1,n)$ and completely contained in the $y$-axis  (reps. $x$-axis) when $i=1$ (resp. $j=n$).
Pedagogical examples are  presented in \pref{fig:flats.hyperplane}.

\begin{proof}[Proof of \pref{lem:flatYieldsEnsemble}]
Recall that an ensemble is the restriction to $\posetz_n$ of  a union
\[
\left(\bigcup_{0\leq i\leq k}  [A_i,B_i] \right)\cap \posetz_n
\]
of intervals $[A_i,B_i]$ of $\poseti_n$ satisfying the following four conditions:
\begin{enumerate}[(1)]
\item $A_0=0;$
\item $A_i\leq B_i$ for $0 \le i \le k$;
\item $B_i +(1,1) \leq A_{i+1}$ for $0 \le i \le k-1$; and
\item $\ell(B_k)< n$ or $B_k=\infty$.  
\end{enumerate}

Consider a $k$-flat $L$ and denote by $E$ the set of extreme rays of $L$. 
 We will prove that $E$ is an ensemble.  Suppose that
\[L = W^0 \cap \bigcap_{(i,j) \in S} \lambda_{i,j}^\perp.\]
Note that
\[\vec{E} = \bigcap_{(i,j) \in S} \erays(\lambda_{i,j}^\perp \cap W^0).\]
Using the structure of extreme rays of hyperplanes $\lambda_{i,j}^\perp$ and the identity 
$$
[a,b]\cap[c,d]=[a\vee c, b\wedge d],
$$
we see immediately that $E \cup \{(0,0)\}$ can be minimally expressed in the form
\[E \cup \{(0,0)\} = [A_0,B_0] \cup \cdots \cup [A_m,B_m].\]
We claim that this expression exhibits $E$ as an ensemble.
 Condition (1) is satisfied because $\{(0,0)\} \in E \cup \{(0,0)\}.$
Condition (2) follows from the minimality of the expression for $E$.
Note that for all $0 \le p \le m-1,$ there exists $(i,j) \in S$ such that $B_p \le \varphi(i,j)$ and $\varphi(i,j) + (1,1) \le A_{p+1}$.  In particular, we have $B_p + (1,1) \le A_p$ for all $0 \le p \le m-1$, which is Condition (3).
Note that either $B_m = \varphi(i,j)$ for some $(i,j) \in S$ for $B_m = \infty$.  Condition (4) follows from the fact that $\ell(\varphi(i,j)) \le n-1$ for all $1 \le i,j \le n$.
Therefore, $E$ is an ensemble.

It remains to prove that $E$ is a $k$-ensemble.  Note that because $E$ is a union of intervals, the number of distinct levels of elements of $E$ is the same as the length of any maximal chain in $E$.  Suppose that $E$ is an $d$-ensemble and that $C \subseteq E$ is an $d$-chain.  \pref{thm:Characterization.Faces} ensures that there exists an $d$-face $F$ whose set of extreme rays is $\vec{C}$.  Because $F$ is a non-negatively spanned by $\vec{C}$ and $L$ is non-negatively spanned by $\vec{E}$, we obtain that $F \subseteq L$.  It follows that $d \le k$.  Consider a $k$-face $F \subseteq L$, which exists because $L$ is a $k$-flat and every flat is the union of the (finitely many) faces it contains.  We see immediately that $\ind{\erays(F)} \subseteq E$.  \pref{thm:Characterization.Faces} guarantees that $\ind{\erays(F)}$ is a $k$-chain in $\posetz_n$, from which it follows that $k \le d$.  Hence, we have $k = d$, so that $E$ is a $k$-ensemble. The lemma follows.
\end{proof}

The proof of \pref{lem:ensembleYieldsFlat} will rely on the following special case.

\begin{lemma}
\label{lem:getSpecialFlats}
Given two points of $A$ and $B$ of  $\poseti_n$ such that $B+(1,1)\leq A$, there exists a flat $L$ of $\glngeo$ such that
\[\left([0, B]\cup [A, \infty]\right) \cap \posetz_n = \erays(L).\]
\end{lemma}

\begin{figure}
   \begin{multicols}{3}

\begin{minipage}{\linewidth}
\centering

                 \includegraphics[scale=1]{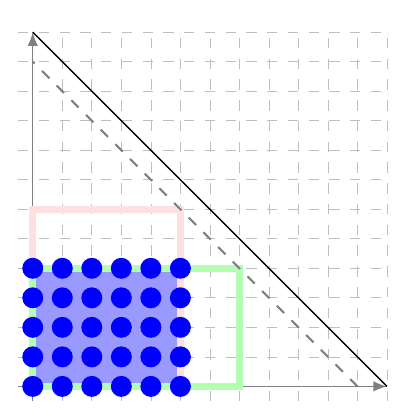}
                 \end{minipage}
\subcaption*{  {\footnotesize (1)  $\lambda^\perp_{c,c}\cap \lambda^\perp_{d,d}$} \label{fig:fig5a}
   }
   
\begin{minipage}{\linewidth}
\centering

                 \includegraphics[scale=1]{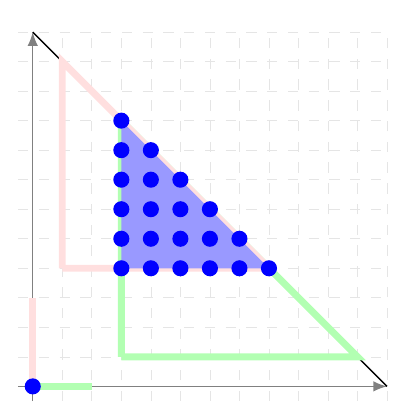}
                 \end{minipage}

              \subcaption*{ {  \footnotesize (2.a)   \quad 
    $\lambda^\perp_{1,i}\cap \lambda^\perp_{j,n}$}  \label{fig:fig5b}
    }
    
    \begin{minipage}{\linewidth}
\centering

                 \includegraphics[scale=1]{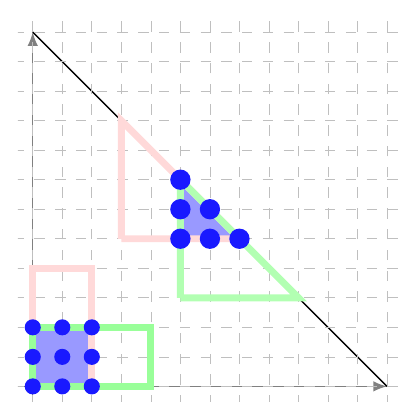}
                 \end{minipage}
\subcaption*{ { \footnotesize (2.b) 
    $\lambda^\perp_{a,b}\cap \lambda^\perp_{c,d}$  ($a>c$, $b<d$) } \label{fig:fig5c}
  }
   \end{multicols}
\caption{
 {\bf  Decomposition of an ensemble  $[(0,0),B] \cup [A,\infty]\cap \posetz_n$} \footnotesize
 with  $A,B$ in $\poseti_n$,  $\ell(B) < n$, and $B+ (1,1) \le A$. The three figures illustrate how to construct a flat $L$ whose set of extreme rays $\erays(L)$ corresponds to the ensemble  $[(0,0),B] \cup [A,\infty]\cap \posetz_n$. 
 In every case, $L$ is  the intersection of two hyperplanes. 
 There are three cases to consider: (Case  1)  $A=\infty$  (Case 2.a) \ $A=\infty$ and $B\neq 0$, (Case 2.b)   $A\neq \infty$ and $B=0$. 
}
\label{fig:BasisForEnsembles}
\end{figure}

The proof of \ref{lem:getSpecialFlats} also relies on the explicit descriptions of the set of extreme rays of a hyperplane of $\glngeo$.

\begin{proof}
The proof is inspired by      \pref{fig:BasisForEnsembles}. Write $B=(x,y) \in \poset_n$.
We divide into cases based on whether $A = \infty$.
\begin{itemize}
\item Case 1: $A=\infty$.  Consider the flat $L = W \cap \lambda^\perp_{x+1,x+1} \cap \lambda^\perp_{n-y,n-y}$ of $\glngeo$.
A direct calculation shows that
\[\left([0, B] \cup [A, \infty]\right) \cap \posetz_n = \erays\left(L\right).\]
\item Case 2: $A \not= \infty$. Write $A=(u, w)\in \poset_n$.
Consider the flat $L = W \cap \lambda^\perp_{x+1, n-w+1} \cap \lambda^\perp_{u,n-y}$ of $\glngeo$. 
A direct calculation shows that
\[\left([0, B] \cup [A, \infty]\right) \cap \posetz_n = \erays\left(L\right).\] 
\end{itemize}
The cases clearly exhaust all possibilities for $A$ and $B$, and the lemma follows.
\end{proof}

To prove \pref{lem:ensembleYieldsFlat}, we exploit the fact that every ensemble is an intersection of ensembles considered by \pref{lem:getSpecialFlats}.

\begin{proof}[Proof of \pref{lem:ensembleYieldsFlat}]
Consider  an ensemble given by 
\[E=\posetz_n \cap \left([A_0, B_0]\cup [A_1, B_1]\cup \cdots \cup [A_k, B_k]\right).\]
Define the intervals 
\begin{align*}
D_i&= (0, B_i]\cup [A_{i+1}, \infty) \text{ for } 0 \le i \le k-1\\
D_k&= (0, B_k].
\end{align*}
\pref{lem:getSpecialFlats} guarantees that each $D_i$ is the set of extreme rays of a flat $F_i$.  Consider the flat
\[F = \bigcap_{i=0}^k F_i.\]
A direct calculation shows that
\[\erays(F) = \bigcap_{i=0}^k \erays(F_i) = E,\]
as desired.
\end{proof}

\pref{thm:Characterization.Flats} follows without much work from Lemmas~\ref{lem:flatYieldsEnsemble} and~\ref{lem:ensembleYieldsFlat}.

\begin{proof}[Proof of \pref{thm:Characterization.Flats}]
\pref{lem:flatYieldsEnsemble} ensures that $\erays$ induces a function from the set of $k$-flats of $\glngeo$ to the set of $k$-ensembles in $\posetz_n$.  Because every flat is non-negatively spanned by its extreme rays, the induced function is injective.

It remains to prove that the induced function is surjective.  Let $E$ be a $k$-ensemble.  \pref{lem:ensembleYieldsFlat} guarantees the existence of a flat $F$ whose set of extreme rays is $E$.  The fact that $F$ is a $k$-flat follows from \pref{lem:flatYieldsEnsemble} and the fact that $E$ is a $k$-ensemble.  Therefore, the induced function is surjective, which implies the theorem.
\end{proof}
  
\subsection{Proof of \pref{thm:flatCount}}
\label{sec:genFuncFlats}

In order to count $k$-ensembles, we generalize ensembles to pseudo-ensembles and relate the counts of ensembles and pseudo-ensembles to one another.

\begin{definition}
For $n \ge 0$, a \emph{pseudo-ensemble} is the restriction to $\poset_n$ of  a union
\[
\left(\bigcup_{0\leq i\leq k}  [A_i,B_i] \right)\cap \poset_n
\]
of intervals $[A_i,B_i]$ of $\poseti_n$ satisfying the following four conditions:
\begin{enumerate}[(1')]
\item $A_i\leq B_i$ for $0 \le i \le k$;
\item $B_i +(1,1) \leq A_{i+1}$ for $0 \le i \le k-1$; and
\item $\ell(B_k)< n$ or $B_k=\infty$.  
\end{enumerate}
We say that $E$ is a \emph{$k$-pseudo-ensemble} if $k = \left|\ell(E)\right| - 1$ so that $k-1$ counts the number of distinct levels of elements of $E$.  We call $A_0$ the \emph{starting point} of $E$.
\end{definition}
\begin{example}
A subset $E \subseteq \posetz_n$ is a $k$-ensemble if and only if $E \cup \{0,0\} \subseteq \poset_n$ is a $k$-pseudo-ensemble (necessarily with starting point 0).  In particular, the function defined by $E \mapsto E \cup \{0,0\}$ induces a bijection from the set of $k$-ensembles in $\posetz_n$ to the set of $k$-pseudo-ensembles in $\poset_n$ with starting point $(0,0)$.
For all $n \ge 0,$ the set $\emptyset \subseteq \poset_n$ is the unique $(-1)$-pseudo-ensemble in $\poset_n$.
\end{example}

Denote by $\rho(n,k)$ the number of $k$-pseudo-ensembles in $\poset_n$ for $n \ge 0$.  Define $\rho(n,k) = 0$ for $n < 0$ and note that $\rho(n,k) = 0$ for $k < -1$ and for $k > n$.
The following two lemmas  illustrate the role of $\rho(n,k)$ in counting ensembles: \pref{thm:flatCount} will follow from the recurrence relations presented in the lemmas by routine algebraic manipulations involving generating functions.

\begin{lemma}
\label{lem:rhoNKrecursion}
For all $n \ge 0$ and all $k$, we have
\[\rho(n,k) = \delta_{k,-1} + \sum_{\ell=0}^n h(n-\ell,k).\]
\end{lemma}

\begin{lemma}
\label{lem:hNKrecursion}
For all $n \ge 0$ and all $k$, we have
\[h(n,k) = \delta_{n,k} + \sum_{\ell=0}^{n-1} (\ell+1)\rho(n-\ell-2,k-\ell-1).\]
\end{lemma}

The idea of \pref{lem:rhoNKrecursion} is to divide into cases based on the starting point of a pseudo-ensemble.  After translating the starting point to $(0,0)$ and removing $(0,0)$, we obtain an ensemble.

\begin{proof}[Proof of \pref{lem:rhoNKrecursion}]
The lemma is clear for $k < 0$.  Assume for the remainder of the proof that $k \ge 0$.  Denote by $A_P$ the set of $k$-pseudo-ensembles in $\poset_n$ with starting point $P$.
Recall that translation by $-P$ induces an isomorphism from $(P,\infty)$ to $\posetz_{n-\ell(P)}$.  We see immediately that the composite of translation by $-P$ and the removal of $(0,0)$ (i.e., the function $S \mapsto (-P + S) \setminus \{(0,0)\}$) induces a bijection from $A_P$ to the set of $k$-ensembles in $\posetz_{n-\ell(P)}$. Hence, we have $\left|A_P\right| = h(n-\ell(P),k).$

Recall that $\posetz_n$ contains $\left|L_\ell\right| = \ell+1$ elements of level $\ell$ for $1 \le \ell \le n$.  Grouping by level, we obtain that
\[\rho(n,k) = \sum_{\ell=0}^n h(n-\ell,k)\]
as desired.
We need to introduce the correction term $\delta_{k,-1}$ in order to deal with the case of $k < 0$.
\end{proof}

The proof of \pref{lem:hNKrecursion} is similar.  We instead divide into cases based on the level of $B_0 \in \poset_n$, the maximum of the first interval in the definition of an ensemble.  Translating $B_0$ to $(-1,-1)$ and intersecting with $\poset_n$, we obtain a pseudo-ensemble.

\begin{proof}[Proof of \pref{lem:hNKrecursion}]
The lemma is clear for $k < 0 $ and for $k \ge n$.  Assume for the remainder of the proof that $0 \le k \le n-1$.
Denote by $B_P$ the set of $k$-ensembles in $\posetz_n$ whose first interval is $[0,P]$.  The hypothesis that $k < n$ ensures that $B_P = \emptyset$, while the definition of an ensemble guarantees that $B_P = \emptyset$ for $\ell(P) = n$.
Recall that translation by $-P$ induces an isomorphism from $(P,\infty)$ to $\posetz_{n-\ell(P)}$.

Recall that translation by $-P-(1,1)$ induces an isomorphism from $[P+(1,1),\infty)$ to $\poset_{n-\ell(P)-2}$.
We see immediately that the composite of translation by $-P-(1,1)$ and intersection with $\mathbb{N}^2$ (i.e., the function $S \mapsto (-P + S) \cap \mathbb{N}^2$) induces a bijection from $C_P$ to the set of $(k-\ell(P)-1)$-pseudo-ensembles in $\posetz_{n-\ell(P)-2}$. Hence, we have $\left|A_P\right| = \rho(n-\ell(P)-2,k-\ell(P)-1).$

Grouping by level, we obtain that
\[h(n,k) = \sum_{\ell=0}^{n-1} (\ell+1)\rho(n-\ell-2,k-\ell-1).\]
as desired.
We need to introduce the correction term $\delta_{k,n}$ in order to deal with the case of $k \notin [0,n-1]$.
\end{proof}

We are now ready to solve the recurrence relations of Lemmas~\ref{lem:rhoNKrecursion} and~\ref{lem:hNKrecursion} in order to prove \pref{thm:flatCount}.

\begin{proof}[Proof of \pref{thm:flatCount}]
Lemmas~\ref{lem:rhoNKrecursion} and~\ref{lem:hNKrecursion} yield the linear equations
\begin{align*}
P(s,t)&=\frac{1}{t(1-s)}+\frac{1}{(1-s)^2}  H(s,t)
\\
H(s,t) &=\frac{1}{(1-st)^2}+\frac{s}{1-st}+\frac{s^2t}{(1-st)^2} P(s,t),
\end{align*}
respectively. Solving for $H(s,t)$ yields the theorem.
\end{proof}

The proof of \pref{cor:flatRecurse} is similar to that of \pref{cor:faceRecurse}.  In particular, \pref{cor:flatRecurse} can be derived by noting that the coefficient of $s^nt^k$ in 
\[(1-2 s+s^2-2 s t+3 s^2 t-2 s^3 t+s^2 t^2-2 s^3 t^2+s^4 t^2) H(s,t)\]
vanishes for $n \ge 4$ due to \pref{thm:flatCount}.


\begin{table}[bht]
\begin{tabular}{llp{12cm}}
& & \\
& & \\
 $|A|$ &:& the cardinality of the set $A$\\
$2^A$ &:& the power set of a set $A$\\
$\binom{n}{k}$ &:& the binomial coefficient  $n$ chosen $k$\\
$[n]$ &:& the set $\{1,2,\cdots, n\}$ \\
$\chi_S(n)$ &:& characteristic function of the set $S$ evaluated at $n$\\
$\mathbb{N}$ &:& the set of non-negative integers\\
$\mathbb{R}$ &:& the set of real numbers \\
$\mathbb{R}_{\geq 0}$ &:& the set of non-negative real numbers \\
$v_{\geq}$ &:& the result of rearranging the components of $v \in \mathbb{R}^n$  in non-increasing order\\
$\mathfrak{h}$ &:& the Cartan sub-algebra of $\mathfrak{gl}_n$ consisting of diagonal matrices \\
$W^0$ &:& the open Weyl chamber of $\mathfrak{gl}_n$\\
$W$ &:& the closure of $W_0$ in $\mathfrak{h}$ \\
$\lambda^\perp$ &:& the hyperplane in $\mathfrak{h}$ on which $\lambda \in \mathfrak{h}^*$ vanishes\\
$V$ &:& the first fundamental representation of $\mathfrak{gl}_n$  (the ``standard representation")\\
$\Exterior^2$ &:& the second  fundamental representation of $\mathfrak{gl}_n$ (the second exterior power of the standard representation)\\
$\lambda_{i,j}$  $(i\leq j)$ &:& the weight $x_i+x_j \in \mathfrak{h}^*$ (resp. $x_i \in \mathfrak{h}^*$) if $i \not= j$ (resp. $i=j$)\\
\\
$\ell(p)$ &:& the level $a+b$ of $p=(a,b) \in \mathbb{N}^2$ \\
$\posetz_n$ &:&  the set of points of level in $[1,n]$ in $\mathbb{N}^2$\\
$\poset_n$ &:&  the set of points of level in $[0,n]$ in $\mathbb{N}^2$\\
$\poseti_n$ &:& the completion $\poset_n\sqcup \{\infty\}$ of $\poset_n$ by a maximum\\
$\pi^S_\ell$ &: & the number of elements of $S\subseteq [n]$ that are at least $n-\ell+1$\\
$\ray{a,b}$ &:& the vector $(\underbrace{1,\ldots ,1}_{a}, \underbrace{0,\ldots ,0}_{n-a-b} , \underbrace{-1,\ldots ,-1}_{b}) \in \mathfrak{h}$\\
$e^s_\ell$ & :& the extreme ray $\ray{\pi_\ell^S,\ell-\pi_\ell^S}$ of the chamber $C(s)$ whose index has level $\ell$, where $S \subseteq [n]$ has characteristic vector $s$\\
\\
$g(n,k)$  & : &  the number of $k$-faces of $I(\mathfrak{gl}_n,V\oplus \Exterior^2)$ \\
$h(n,k)$  & : & the number of $k$-flats of $I(\mathfrak{gl}_n,V\oplus \Exterior^2)$ \\
$G(s,t)$  & : & the generating function $\sum_{n,k\ge0} g(n,k)s^nt^k$ \\
$H(s,t)$  & : & the generating function $\sum_{n,k\ge0} h(n,k)s^nt^k$ \\ 

\end{tabular}
\vspace{.5 cm}
\caption{List of symbols.}
\label{tab:notations}
\end{table}

\end{document}